\documentclass[12pt]{article}
\usepackage{amsmath}
\def\blacksquare{{\vcenter{\hbox{\vrule width1.3ex height1.3ex}}}}
\def\qed{\hfill$\blacksquare$}
\def\tr{{\bf tr}\,}
\def\re{{\bf Re}\,}
\def\cU{{\mathcal U}}

\def\ok{{\oplus 0_{n-2}}}
\def\oo{{\oplus 0}}
\begin{document}
\newtheorem{theorem}{Theorem}
\newtheorem{lemma}[theorem]{Lemma}
\newtheorem{corollary}[theorem]{Corollary}
\newtheorem{example}{Example}
\title{New pairs of matrices with convex generalized numerical ranges}
\author{Wai-Shun Cheung}
\date{}
\maketitle

\begin{abstract}

In this article, we are going to search for $n\times n$ matrices $A$ and $B$ such that their generalized numerical range
$$W_A(B)=\{\tr AU^*BU \ :\ U^*U=UU^*=I\}$$
is convex.  More specifically, we consider $A=\hat{A}\ok$ and $B=\hat{B}\ok$ where $\hat{A}$ and $\hat{B}$ are $2\times 2$. If $W_A(B)=W_{\hat{A}}(\hat{B})$ then it is a convex set. 

\end{abstract}



\section{Introduction}
Let $M_n$ be the space of all $n\times n$ matrices with standard basis $\{E_{11}, E_{12}, \ldots, E_{nn}\}$, and $U_n$ be the group of all $n\times n$ unitary matrices.  

For $B\in M_n$, the classical numerical range of $B$ is the set 
$$W(B)=\{x^*Bx\ :\  x\mbox{ is a unit vector}\}.$$
The classical numerical range is a compact set which contains all the eigenvalues of $B$, and it is a convex set by the famous Toeplitz-Hausdorff Theorem \cite{H, To}. See \cite[Chapter 1]{HJ2} for a nice discussion.

Note that $W(B)=\{\tr E_{11}X\ :\ X\in U(B)\}$ where $U(B)=\{VBV^*\ :\ V\in U_n\}$ is the unitary orbit of $B$.   This inspires the following generalization. Let $C\in M_n$, the set
$$W_C(B)=\{\tr CX\ :\ X\in U(B)\}$$
is called the $C$-numerical range of $B$.  Therefore the classical numerical range of $B$ is the $E_{11}$-numerical range of $B$.  Note that the $C$-numerical range of $B$ is the $B$-numerical range of $C$.   

In 1975, Westwick \cite{W} showed that if $C$ is Hermitian then $W_C(B)$ is convex. (See another proof by Poon \cite{P}.) Hence $W_C(B)$ is convex if $C$ is a normal matrix with collinear eigenvalues. Conjectured by Marcus \cite{M} in 1975 and confirmed by Au-Yeung and Tsing \cite{AT} in 1983, if $C$ is normal and $W_C(B)$ is convex for all $B$ then the eigenvalues of $C$ must be collinear.

In 1984, Tsing \cite{T1} proved that if $C$ is rank one then $W_C(B)$ is convex for all $B$.  A consequence is that $W_C(B)$ is convex for any $B,C\in M_2$.

\medskip
Problem 1. \it Find more $B,C$ with convex $W_C(B)$. \rm

\medskip
Problem 2. \it So far, for all $B,C$ with convex $W_C(B)$, one of $B$ and $C$ must have collinear eigenvalues. Is it a general rule? \rm

\medskip
In 1991, Li and Tsing \cite{LT} showed that if $C=\lambda I+C_0$ where $C_0$ is the block-shift form matrix then $W_C(B)$ is always a circular disc centered at $\lambda\tr B$. Indeed, if $W_C(C^*)$ is a circular disc centered at $0$ then $C$ must be a shift-block form matrix. 

\medskip
Problem 3. \it Suppose $W_C(C^*)$ is a circular disc. Does $C=\lambda I+C_0$ where $C_0$ is a block-shift from matrix? \rm

\medskip
Although $W_C(B)$ may fail to be convex, it is proved in 1981 by Tsing \cite{T} that if $C$ is normal then $W_C(B)$ is star-shaped. Later in 1996, Cheung and Tsing \cite{CT} showed that $W_C(B)$ is star-shaped for all $C$ and $B$.

In this article, we will do Problem 1, i.e., to search other pairs of $B,C\in M_n$ such that $W_C(B)$ is convex.  More specifically, we consider:

\medskip
Problem 4. \it Find $A,B\in M_2$ such that $W_{A\oplus 0_k}(B\oplus 0_k)=W_A(B)$. \rm

\medskip
If we can find such $A,B$ then $(A\oplus 0_k, B\oplus 0_k)$ is a ``convex pair''.

In the end, we will answer Problem 2 and Problem 3 as well.

Let's have some more notations:  For $B=\begin{pmatrix}b_{11}&b_{12}\\b_{21}&b_{22}\end{pmatrix}\in M_2$, we write
$$B_0=B-\left(\frac12 \tr B\right)I \quad\hbox{and}\quad
B(\epsilon)=
\begin{pmatrix}1&0\\0&\epsilon\end{pmatrix}B
\begin{pmatrix}1&0\\0&\epsilon\end{pmatrix}
=\begin{pmatrix}b_{11}&\epsilon b_{12}\\\epsilon b_{21}&\epsilon^2b_{22}\end{pmatrix}$$
where $0\le\epsilon\le 1$.

Recall that the numerical radius of a square matrix $A$ is given by
$$r(A)=\max_{x\in W(A)}|x|.$$
%  We also define
%$$\cS_n=\{(A,B)\in M_n\times M_n\ :\ W_A(B)\mbox{ is convex}\}.$$

\section{Lemmas}
In 1932, Murnaghan \cite{Mu} proved the original Elliptical Range Theorem, which states that the classical numerical range of $A\in M_2$ is an elliptical disc centered at $\dfrac12\tr A$ and the two eigenvalues are the foci on the major axes.  In 1994, Nakasato \cite{N} generalized it to general $W_C(A)$. Let's state Nakasato's result as our first lemma. 

\begin{lemma} \label{L0} {\rm (Elliptical Range Theorem)} Let $A,B\in M_2$. If $A=\mu\begin{pmatrix}a&a_{12}\\a_{21}&a\end{pmatrix}$ and $B=\nu\begin{pmatrix}b&b_{12}\\b_{21}&b\end{pmatrix}$ with $a_{12}\ge a_{21}\ge 0$ and $b_{12}\ge b_{21}\ge 0$ then 
$$W_A(B)=2\mu\nu W\left(\begin{pmatrix}ab&a_{12}b_{12}\\a_{21}b_{21}&b\end{pmatrix}\right)=2ab+\cU$$
which is an elliptical disc centered at $2ab$ and
$$\cU=W_A(B_0)=W_B(A_0)=W_{A_0}(B_0)$$
is an ellipitical disc centered at $0$.
\end{lemma}

Please also see another proof of Lemma~\ref{L0} by  Li \cite{Li}.

\begin{lemma} \label{L1} {\rm \cite{CT}} Let $A,B=\begin{pmatrix}b_{11}&b_{12}\\b_{21}&b_{22}\end{pmatrix}\in M_2$. We have $W_A\left(\begin{pmatrix}b_{11}&\epsilon b_{12}\\\epsilon b_{21}&b_{22}\end{pmatrix}\right)\subseteq W_A(B)$ where $0\le\epsilon\le 1$.
\end{lemma}

From now on, we always assume that $n\ge 3$.

\begin{lemma} \label{L2} Let $A,B\in M_2$ and $U\in U_n$. There exists $\hat{A}\in U(A)$ and $\hat{B}\in U(B)$ such that 
$$\tr (A\ok)U^*(B\ok)U=\alpha\tr \hat{A}(\epsilon)\hat{B}=\alpha\tr \hat{A}\hat{B}(\epsilon)$$
for some $0\le \alpha,\epsilon\le 1$.  If $n=3$, then $\alpha=1$.

Consequently, we have
$$W_{A\ok}(B\ok)=\bigcup_{\hat{A}\in U(A), 0\le\alpha,\epsilon\le 1}\alpha W_B(\hat{A}(\epsilon)))=\bigcup_{\hat{B}\in U(B), 0\le\alpha,\epsilon\le 1}\alpha W_A(\hat{B}(\epsilon))$$ when $n\ge 4$, and
$$W_{A\oplus 0}(B\oplus 0)=\bigcup_{\hat{A}\in U(A), 0\le\epsilon\le 1}W_B(\hat{A}(\epsilon)))=\bigcup_{\hat{B}\in U(B), 0\le\epsilon\le 1}W_A(\hat{B}(\epsilon)).$$

\end{lemma}

{\it Proof.} It follows from the singular value decomposition of the leading $2\times 2$ prinicpal submatrix of $U$. \qed

\begin{lemma} \label{L5} Let $A,B\in M_2$.  If $W_{A\ok}(B\ok)=W_A(B)$ then $0\in W_A(B)$.  Consequently if $W_{A\ok}(B\ok)=W_A(B)$ for some $n$, then it is true for all $n$.
\end{lemma}

{\it Proof.} Without loss of generality, we assume $A=\begin{pmatrix}a_{11}&a_{12}\\a_{21}&a_{11}\end{pmatrix}$ where $a_{11}=\frac12\tr A$.  By  Lemma~\ref{L2}, we have
\begin{eqnarray*}
&&\left(\frac{1+\epsilon^2}{2}\right)a_{11}\tr B+W_{B-\left(\frac12\tr B\right)I}(A(\epsilon))=W_B(A(\epsilon))\\
&\subseteq& W_B(A)=a_{11}\tr B+W_{B-\left(\frac12\tr B\right)I}(A)
\end{eqnarray*}
hence
\begin{equation} \label{eqn1}
\left(\frac{-1+\epsilon^2}{2}\right)a_{11}\tr B+W_{B-\left(\frac12\tr B\right)I}(A(\epsilon))
\subseteq W_{B-\left(\frac12\tr B\right)I}(A).
\end{equation}
Note that $$W_{B-\left(\frac12\tr B\right)I}(A)=W_{B-\left(\frac12\tr B\right)I}\left(\begin{pmatrix}0&a_{12}\\ a_{21}&0\end{pmatrix}\right).$$
Note also that
$A(\epsilon)=\begin{pmatrix}a_{11}&\epsilon a_{12}\\\epsilon a_{21}&\epsilon^2 a_{11}\end{pmatrix}$ is unitarily similar to $\begin{pmatrix}\epsilon^2 a_{11}&\epsilon a_{12}\\\epsilon a_{21}&a_{11}\end{pmatrix}$, and so
\begin{eqnarray*}
\epsilon\, W_{B-\left(\frac12\tr B\right)I}(A)
&=&W_{B-\left(\frac12\tr B\right)I}\left(\begin{pmatrix}0&\epsilon a_{12}\\\epsilon a_{21}&0\end{pmatrix}\right)\\
&=&W_{B-\left(\frac12\tr B\right)I}\left(\begin{pmatrix}\frac12(1+\epsilon^2) a_{11}&\epsilon a_{12}\\\epsilon a_{21}&\frac12(1+\epsilon^2) a_{11}\end{pmatrix}\right)\\
&\subseteq& W_{B-\left(\frac12\tr B\right)I}(A(\epsilon)).
\end{eqnarray*}
Therefore (\ref{eqn1}) is possible only if 
$$\left(\frac{-1+\epsilon^2}{2}\right)a_{11}\tr B \in (1-\epsilon)W_{B-\left(\frac12\tr B\right)I}(A)$$
which implies
$$\left(\frac{-1-\epsilon}{2}\right)a_{11}\tr B \in W_{B-\left(\frac12\tr B\right)I}(A)=W_B(A)-a_{11}\tr B.$$
Setting $\epsilon=1$, we have $0\in W_B(A)$. \qed

\begin{lemma} \label{L3}  Let $A,B\in M_2$. If $W_{A\oplus 0_{n-2}}(B\oplus 0_{n-2})=W_A(B)$ then  the largest possible $\alpha$ such that $\alpha W(A)W(B)\subseteq W_A(B)$ satisfies $1\le \alpha\le 4$.
\end{lemma}

{\it Proof.} Let $a\in W(A)$, then there exists $\hat{A}=\begin{pmatrix}a&a_{12}\\a_{21}&a_{22}\end{pmatrix}\in U(A)$.  By Lemma~\ref{L2}, we have $aW(B)=W_{\hat{A}(0)}(B)\subseteq W_A(B)$.

By Lemma~\ref{L0}, we have $W_A(B)=2W(\hat{A}\circ \hat{B})$ for some $\hat{A}\in U(A)$ and $\hat{B}\in U(B)$. Therefore, if $\alpha W(A)W(B)\subseteq W_A(B)$, we have $\alpha r(A)r(B)\le 2 r(\hat{A}\circ \hat{B})$.  However, $r(\hat{A}\circ\hat{B})\le 2r(\hat{A})r(\hat{B})$ \cite[Corollary 1.7.25]{HJ2}. Therefore, $\alpha r(A)r(B) \le 4r(A)r(B)$ and thus $\alpha\le 4$. 
\qed

\medskip
The lower bound and the upper bound for $\alpha$ are both sharp. The lower bound is sharp because $W_{E_{11}}(E_{11})=[0,1]=W(E_{11})W(E_{11})$ and the upper bound is sharp because $W_{E_{12}}(E_{12})=4W(E_{12})W(E_{12})$.

\begin{lemma} \label{L4} Let $A,B\in M_2$. If $W_{A\oplus 0_{n-2}}(B\oplus 0_{n-2})=W_A(B)$ then
$$W_B(A_0)=W_{B\ok}(A_0\ok)=W_A(B_0)=W_{A_0}(B_0\ok).$$
\end{lemma}

{\it Proof.} By Lemma~\ref{L2}, we have $W_{\hat{A}(\epsilon)}(B)\subseteq W_A(B)$ for any $\hat{A}\in U(A)$, which implies, by Lemma~\ref{L0}, that 
$$W_{\hat{A}(\epsilon)}\left(B-\left(\frac12\tr B\right)I\right)\subseteq W_A\left(B-\left(\frac12\tr B\right)I\right)$$
 for any $\hat{A}\in U(A)$, and then by Lemma~\ref{L2} again, we have  
$$W_A\left(B-\left(\frac12\tr B\right)I\right)=W_{A\oplus 0_{n-2}}\left(\left(B-\left(\frac12\tr B\right)I\right)\oplus 0_{n-2}\right).$$
\qed

\section{Main Results}
An implication of Lemma~\ref{L5} is that we only need to consider the case $n=3$.  

First of all, we have a sufficient condition.

\begin{theorem} \label{M0}
Let $A,B\in M_2$. Suppose 
$$(\tr A)W(B)+(\tr B)W(A)-(\tr A)(\tr B)\subseteq W_{A}(B),$$
or equivalently
$$(\tr A)W(B_0)+(\tr B)W(A_0)-\frac12(\tr A)(\tr B)\subseteq W_{A_0}(B_0)$$
then  $W_A(B)=W_{A\oo}(B\oo)$.
\end{theorem}

{\it Proof.} Let $\hat{A}=\begin{pmatrix}a_{11}&a_{12}\\a_{21}&\tr A-a_{11}\end{pmatrix}\in U(A)$ and 
$\hat{B}=\begin{pmatrix}b_{11}&b_{12}\\b_{21}&\tr B-b_{11}\end{pmatrix}\in U(B)$. We have
\begin{eqnarray*}
&&\tr(\hat{A}\hat{B}(\epsilon))\\
&=& a_{11}b_{11}+\epsilon (a_{12}b_{21}+a_{21}b_{12})+\epsilon^2(\tr A-a_{11})(\tr B-b_{11})\\
&=&\frac12(1+\epsilon^2)\left(a_{11}b_{11}+\frac{2\epsilon}{1+\epsilon^2}(a_{12}b_{21}+a_{21}b_{12})+(\tr A-a_{11})(\tr B-b_{11})\right)\\
&&+\frac12(1-\epsilon^2)(a_{11}\tr B+b_{11}\tr A-(\tr A)(\tr B))\\
&=&\frac12(1+\epsilon^2)\tr\left(\hat{A}\begin{pmatrix}b_{11}&\frac{2\epsilon}{1+\epsilon^2} b_{12}\\\frac{2\epsilon}{1+\epsilon^2} b_{21}&\tr B-b_{11}\end{pmatrix}\right)\\
&&+\frac12(1-\epsilon^2)(a_{11}\tr B+b_{11}\tr A-(\tr A)(\tr B))\\
&\in& W_A(B) \quad\mbox{\rm by Lemma~\ref{L1} and the assumption}.
\end{eqnarray*}
Hence by Lemma~\ref{L2}, we have $W_A(B)=W_{A\oo}(B\oo)$. \qed

\medskip
Let's replace the condition in Theorem~\ref{M0} with a stronger one to make it easier to apply.

\begin{theorem} \label{M4}
Let $A,B\in M_2$. Suppose the disc centered at $0$ with radius
$|\tr A|r(B_0)+|\tr B|r(A_0)+\frac12|\tr A||\tr B|$ lies inside $W_{A_0}(B_0)$, then $W_A(B)=W_{A\oo}(B\oo)$.
\end{theorem}

{\it Proof.} It is a direct consequence of Theorem~\ref{M0}.
\qed

\medskip
It turns out that the sufficient condition is also necessary if one of the two matrices is trace $0$.

\begin{theorem} \label{M1} Let $A,B\in M_2$. Suppose $\tr B=0$. Then $W_A(B)=W_{A\oo}(B\oo)$ iff $(\tr A)r(B)\subseteq W_A(B)$.
\end{theorem}

{\it Proof.}  Sufficency follows from Theorem~\ref{M0}.

Suppose that $W_A(B)=W_{A\oo}(B\oo)$.
If there exists $\beta\in W(B)$ such that $(\tr A)\beta \notin W_A(B)$, then there exists $\theta$ such that $\re e^{i\theta}\tr(A)\beta>\re e^{i\theta}x$ for all $x\in W_A(B)$. Without loss of generality, we can assume that $\hat{A}=\begin{pmatrix}a_{11}&a_{12}\\a_{21}&a_{22}\end{pmatrix}$ and $\hat{B}=\begin{pmatrix}\beta&b_{12}\\b_{21}&-\beta\end{pmatrix}$ are such that $\re e^{i\theta}\tr (\hat{A}\hat{B})$ is the largest possible in $W_A(B)$. 

As $W_A(B)=W_{A\oo}(B\oo)$, we know that 
$\tr (\hat{A}\hat{B}(\epsilon))\in W_A(B)$
and hence
$$\re e^{i\theta}(a_{11}\beta+\epsilon(a_{12}b_{21}+a_{21}b_{12})-\epsilon^2a_{22}\beta)
\le \re e^{i\theta}(a_{11}\beta+(a_{12}b_{21}+a_{21}b_{12})-a_{22}\beta).$$
Reorganizing, we have
$$\re e^{i\theta}((1-\epsilon)(a_{12}b_{21}+a_{21}b_{12})-(1-\epsilon^2)a_{22}\beta)\ge 0$$
and so
$$\re e^{i\theta}((a_{12}b_{21}+a_{21}b_{12})-(1+\epsilon)a_{22}\beta)\ge 0.$$
Thus, put $\epsilon=1$, we have
$$\re e^{i\theta}((a_{12}b_{21}+a_{21}b_{12})-a_{22}\beta)\ge \re e^{i\theta}a_{22}\beta.$$
Therefore
$$\re e^{i\theta} \tr (\hat{A}\hat{B})\ge \re e^{i\theta}(a_{11}\beta+a_{22}\beta)=\re e^{i\theta}(\tr A)\beta$$
which is a contradiction. \qed

\medskip
As a corollary, we have a necessary condition.

\begin{theorem} \label{M2} Let $A, B\in M_2$. If $W_A(B)=W_{A\oo}(B\oo)$ then 
$$(\tr A)W(B)\cup(\tr B)W(A)\subseteq W_A(B).$$
\end{theorem}

{\it Proof.} By Lemma~\ref{L4}, we know 
$$W_A\left(B-\left(\frac12\tr B\right)I\right)=W_{A\oplus 0_{n-2}}\left(\left(B-\left(\frac12\tr B\right)I\right)\oplus 0_{n-2}\right),$$  
and then by Theorem~\ref{M1}, we have
\begin{eqnarray*}
&&(\tr A)W(B)-\dfrac12(\tr A)(\tr B)=(\tr A)W\left(B-\left(\frac12\tr B\right)I\right)\\
&\subseteq& W_A\left(B-\left(\frac12\tr B\right)I\right)=W_A(B)-\dfrac12(\tr A)(\tr B).
\end{eqnarray*}
Likewise we have $(\tr B)W(A)\subseteq W_A(B)$. \qed

\medskip
If one of the matrices is Hermitian, then we also have a necessary and sufficient condition.

\begin{theorem} \label{M3} Let $A, B\in M_2$. Suppose $B$ is Hermitian. Then $W_A(B)=W_{A\oo}(B\oo)$ is equivalent to $(\tr A)W(B)\cup(\tr B)W(A)\subseteq W_A(B)$. If, in addition,  both $A$ and $B$ are nonzero matrices, then it is also equivalent to $0\in W(A)\cap W(B)$.
\end{theorem}

{\it Proof.}  It suffices to consider $A,B\ne 0$. Let $b_1\ge b_2$ be the two eigenvalues of $B$ and then $W_A(B)=(b_2-b_1)W(A)+b_1\tr A=(b_1-b_2)W(A)+b_2\tr A$.
By Theorem~\ref{M2},  $W_A(B)=W_{A\oo}(B\oo)$ implies $(\tr A)W(B)\cup(\tr B)W(A)\subseteq W_A(B)$.

Suppose $(\tr A)W(B)\cup(\tr B)W(A)\subseteq W_A(B)$.  Therefore $(b_1+b_2)W(A)=(\tr B)W(A)\subseteq W_A(B)$ which, by considering the area of the two sets, implies that $|b_1+b_2|\le b_1-b_2$ and hence $b_1\ge 0\ge b_2$. $b_1\tr A\in W_A(B)=(b_2-b_1)W(A)+b_1\tr A$ imples that $0\in W(A)$.

$0\in W(A)\cap W(B)$ implies that $b_1\ge 0\ge b_2$ and $0\in W(0.5(e^{i\theta}A+e^{-i\theta}A^*))$ and consequently $W_A(B)=W_{A\oo}(B\oo)$. \qed

\section{New Convex Pairs}

\begin{corollary} Let $A, B\in M_2$ such that $\tr A=\tr B=0$.  We have $W_A(B)=W_{A\oo}(B\oo)$.  Moreover, if both $A$ and $B$ are nonzero, then the largest possible $\alpha$ such that $\alpha W(A)W(B)\subseteq W_A(B)$ satisfies $2\le \alpha\le 4$.
\end{corollary}

{\it Proof.} That $W_A(B)=W_{A\oo}(B\oo)$ follows from Theorem~\ref{M1}.  

Let $a\in W(A)$ and $b\in W(B)$, then there exists $\hat{A}=\begin{pmatrix}a&a_{12}\\a_{21}&-a\end{pmatrix}\in U(A)$ and 
$\hat{B}=\begin{pmatrix}b&b_{12}\\b_{21}&-b\end{pmatrix}\in U(B)$. Hence by Lemma~\ref{L1}, we have
$$2ab=\tr\left(\begin{pmatrix}a&0\\0&-a\end{pmatrix}\begin{pmatrix}b&0\\0&-b\end{pmatrix}\right)\in W_A(B).$$
Hence we have $2W(A)W(B)\subseteq W_A(B)$. 

The upper bound follows from Lemma~\ref{L3}. \qed

\medskip
The lower bound and the upper bound for $\alpha$ are both sharp.  The lower bound is sharp since $W_{E_{11}-E_{22}}(E_{11}-E_{22})=2W(E_{11}-E_{22})W(E_{11}-E_{22})$. The upper bound is sharp, the example is the same as that of Lemma~\ref{L3}.

\begin{example}
$W\left(\begin{pmatrix}2+i&3&0\\1-2i&-2-i&0\\0&0&0\end{pmatrix}, \begin{pmatrix}1+i&2-i&0\\1-2i&-1-i&0\\0&0&0\end{pmatrix}\right)$
is convex.
\end{example}

\begin{corollary} Let $A, B\in M_2$. If $\tr A=0$ and $\frac12\tr B\in W(B_0)$ then $W_A(B)=W_{A\oo}(B\oo)$.
\end{corollary}

{\it Proof.} It follows Theorem~\ref{M1} and that $2W(A)W(B_0)\subseteq W_A(B_0)$. \qed

\begin{example}
$W\left(\begin{pmatrix}2+i&3&0\\1-2i&-2-i&0\\0&0&0\end{pmatrix}, \begin{pmatrix}1&2-i&0\\1-2i&-3&0\\0&0&0\end{pmatrix}\right)$
is convex.
\end{example}

\begin{corollary} If $A=B=\begin{pmatrix}a&b\\0&a\end{pmatrix}$ with 
$|a|\le \dfrac{\sqrt 3-1}{2}|b|$
then $W_A(B)=W_{A\oo}(B\oo)$ is a circular disc centered at $2a^2$ of radius $|b|^2$.
\end{corollary}

{\it Proof.} Without loss of generality, we let $0\le a\le \dfrac{\sqrt 3-1}{2}$ and $b=1$.  Note that $W(A_0)=W(B_0)$ is a circular disc centered at $0$ of radius $0.5$ and $W_{A_0}(B_0)$ is a unit disc centered at $0$.

For $v\in (\tr A)W(B_0)+(\tr B)W(A_0)-\frac12(\tr A)(\tr B)=2a(W(B_0)+W(A_0)-a)$, we have
$$|v|\le 2|a|(2r(A_0)+|a|)=2|a|(1+|a|)\le 1,$$
thus $v\in W_{A_0}(B_0)$.

By Theorem~\ref{M4}, $W_A(B)=W_{A\oo}(B\oo)$. \qed

\begin{example}
If $C=\begin{pmatrix}0.7&2&0\\0&0.7&0\\0&0&0\end{pmatrix}$
then $W(C,C^*)$ is a circular disc centered at $0.98$. This answers \rm Problem 3.
\end{example}

\begin{corollary} If $A=B=\begin{pmatrix}a&b\\1&a\end{pmatrix}$ with 
$b>1$ and $2|a|^2+(1+b)|a|-(b^2-1)\le 0$
then $W_A(B)=W_{A\oo}(B\oo)$
\end{corollary}

{\it Proof.} Note that $r(A_0)=r(B_0)=\frac{1+b}2$.  By Lemma~\ref{L0}, we know that for the circular disc centered at $0$ with radius $b^2-1$ is contained in $W_{A_0}(B_0)$.

 For  $v\in (\tr A)W(B_0)+(\tr B)W(A_0)-\frac12(\tr A)(\tr B)=2a(W(B_0)+W(A_0)-a)$, we have
 $$|v|\le 2|a|(2r(A_0)+|a|)=|a|(1+b+2|a|)\le b^2-1,$$
 thus  $v\in W_{A_0}(B_0)$.

By Theorem~\ref{M4}, $W_A(B)=W_{A\oo}(B\oo)$. \qed

\begin{example}
If $C=\begin{pmatrix}i&3&0\\1&i&0\\0&0&0\end{pmatrix}$
then $W(C,C)$ is convex. Note that the eigenvalues of $C$ are not collinear, this answers \rm Problem 2.
\end{example}

\section*{Acknowledgement}

I wish to thank the referee for his valuable comments.

\medskip
\sc

Department of Mathematics, University of Hong Kong, Hong Kong.

\rm

E-mail addresses: cheungwaishun@gmail.com
\end{document}